\documentclass[12pt]{article}
\usepackage{amsmath,amssymb}
\newtheorem{theorem}{Theorem}[section]
\newtheorem{proposition}[theorem]{Proposition}
\newtheorem{lemma}[theorem]{Lemma}
\def\bull{\vrule height .9ex width .8ex depth -.1ex}
\newenvironment{proof}{\smallbreak \noindent {\bf Proof.~}}
              {\unskip\nobreak\hfill\hskip 2em \bull\par\medbreak}
\newenvironment{proofof}[1]{\medbreak\noindent{\bf Proof of~#1.~}}
              {\unskip\nobreak\hfill\hskip 2em \bull\par\medbreak}
\def\bC{\mathbb{C}}
\def\bN{\mathbb{N}}
\def\bQ{\mathbb{Q}}
\def\bZ{\mathbb{Z}}
\def\cA{\mathcal{A}}
\def\cF{\mathcal{F}}
\def\cT{\mathcal{T}}
\def\cZ{\mathcal{Z}}
\def\cU{\mathcal{U}}
\def\la{\lambda}
\def\si{\sigma}
\def\om{\omega}

\def\CF{\mathop{\mathrm{CF}}}
\def\EP{\mathop{\mathrm{EP}}}
\def\lcm{\mathop{\mathrm{lcm}}}
\title{On a substitution subshift related to\\ the Grigorchuk group}
\author{Yaroslav Vorobets\thanks{%
        Partially supported by the NSF grant DMS-0701298.}}
\date{}
\begin{document}

\maketitle

\section{Introduction}\label{intro}

Let $\tau$ be a substitution over the alphabet $\{a,b,c,d\}$ defined by
relations
$$
\tau(a)=aca, \quad \tau(b)=d, \quad \tau(c)=b, \quad \tau(d)=c.
$$
The substitution
acts on words (finite sequences) over this alphabet as well as on infinite sequences.
It is easy to observe that $\tau$ has a unique invariant sequence $\om$ that is the
limit of words $\tau^k(a)$, $k=1,2,\dots$.  Let $\Omega$ be the smallest
closed set of one-sided infinite sequences over the alphabet $\{a,b,c,d\}$ that contains
$\om$ and is invariant under the shift $\si$ ($\si$ acts on sequences by deleting
the first element).  The set $\Omega$ consists of those sequences for which any finite
subword appears somewhere in $\om$.  The restriction of $\si$ to $\Omega$ is called
a subshift.  Since $\om$ is a fixed point of a substitution, this particular subshift
is called a {\em substitution subshift}.

The substitution $\tau$ plays an important role in the study of the Grigorhuk
group (see the survey \cite{G}).  The Grigorchuk group $G$ is a finitely generated
infinite group where all elements are of finite order.  It has many other remarkable
properties as well.  The group has four generators $a,b,c,d$.  An important fact is that the
substitution $\tau$ gives rise to a homomorphism of $G$ to itself.  It follows that $\tau$
transforms any relator for $G$ into another relator.  Although $G$ has no finite presentation,
it admits a recursive presentation obtained from a finite set of relators by repeatedly
applying $\tau$:
$$
G=\langle a,b,c,d\mid 1=a^2=b^2=c^2=d^2=bcd=\tau^n((ad)^4)=\tau^n((adacac)^4),\, n\ge0\rangle.
$$
The structure of the Grigorchuk group is not completely understood yet.  In view of the above
presentation, it is believed that properties of the sequence $\om$ and the subshift
$\si|_\Omega$ might give an insight on that matter.  In this paper we study dynamics of
$\si|_\Omega$.

\begin{theorem}\label{main1}
The subshift $\si|_\Omega$ is, up to a countable set, continuously conjugated to
the binary odometer.
\end{theorem}

Theorem \ref{main1} suggests that ergodic properties
of the subshift $\si|_\Omega$ are the same as ergodic properties of the binary
odometer.  The following theorem is a detailed version of this suggestion.

\begin{theorem}\label{main2}
(i) The subshift $\si|_\Omega$ has a unique invariant Borel probability measure $\mu$.

(ii) $\si|_\Omega$ is ergodic with respect to the measure $\mu$.  Moreover, every
orbit of $\si|_\Omega$ is uniformly distributed in $\Omega$ with respect to $\mu$.

(iii) $\si|_\Omega$ has purely point spectrum, the eigenvalues being all roots of
unity of order $1,2,\dots,2^k,\dots$.  Each eigenvalue is simple.

(iv) All eigenfunctions of $\si|_\Omega$ are continuous.
\end{theorem}

It turns out that $\om$ is a one-sided analog of what is called Toeplitz
sequences (see, e.g., \cite{D}).  The Toeplitz sequences can be informally described as almost
periodic.  Subshifts generated by Toeplitz sequences are known to be
continuous extensions of odometers.  A nontrivial feature of $\om$ is
that the extension is one-to-one up to a countable set.  In Section \ref{sub}
we place $\om$ into a class of Toeplitz sequences that are as close to periodic
as possible.  Theorems \ref{main1} and \ref{main2} hold for all sequences in that
class.

The paper is organized as follows.  Section \ref{odo} is a survey on odometers.
We discuss their dynamics and determine when a particular odometer is a
continuous factor of a particular topological dynamical system.  In Section
\ref{toe} we collect necessary information about Toeplitz sequences and the
associated subshifts.  In Section \ref{sub}, these results are applied to the
class of Toeplitz sequences that contains $\om$.  The paper ends with the
proof of Theorems \ref{main1} and \ref{main2}.  It should be noted that most
results in Sections \ref{odo} and \ref{toe} are well known to specialists
(although it might not be easy to locate them in the literature).  For reader's
convenience, we include all proofs so that the paper is self-contained.
In the submitted version of the paper the proofs will be replaced by appropriate
references.

\section{Odometers}\label{odo}

In this section we consider general transformations $T:X\to X$ such that
$X$ is a compact topological space and $T$ is a continuous map (not
necessarily one-to-one or onto).  For any $x\in X$ the sequence
$x,Tx,T^2x,\dots$ is called the {\em orbit\/} of the point $x$ under the
transformation $T$.  By $Z(x,T)$ we denote the closure of the orbit.  Note
that $Z(x,T)$ is the smallest closed subset of $X$ that contains $x$ and is
invariant under $T$.  If $C$ is a compact subset of $X$ invariant under
$T$, then $T|_C$ denotes the restriction of $T$ to $C$.  The transformation
$T$ is called {\em transitive\/} if there exists a dense orbit, that is, if
$Z(x,T)=X$ for some $x\in X$.  $T$ is called {\em minimal\/} if each orbit
is dense.

Let $T_1:X_1\to X_1$ and $T_2:X_2\to X_2$ be continuous transformations of
compact sets.  Suppose there exists a continuous map $f:X_1\to X_2$ such
that $f$ is onto and $fT_1=T_2f$ so that the following diagram is
commutative:
$$
\begin{array}{ccc}
X_1 & \stackrel{T_1}{\longrightarrow} & X_1\\[0.2em]
f\Big\downarrow\phantom{f} && \phantom{f}\Big\downarrow f\\[0.25em]
X_2 & \stackrel{T_2}{\longrightarrow} & X_2
\end{array}
$$
Then $T_2$ is called a {\em (continuous)
factor\/} of $T_1$ while $T_1$ is called a {\em (continuous) extension\/}
of $T_2$.  If, in addition, we can choose $f$ to be a homeomorphism then
$T_1$ and $T_2$ are called {\em (continuously) conjugated\/} and $f$ is
called a {\em conjugacy}.

A transformation $T:X\to X$ is called a {\em cyclic permutation\/} if $X$
is a finite set with the discrete topology, $T$ is one-to-one, and an orbit
of $T$ contains all elements of $X$.  Let $n$ denote the cardinality of
$X$.  Then for any $x\in X$ the sequence $x,Tx,T^2x,\dots,T^{n-1}x$ is a
complete list of elements of $X$.  Besides, $T^nx=x$ so that $T$ has order
$n$.  Any cyclic permutation is determined by its order up to conjugacy.

It is easy to observe that a cyclic permutation $T_1$ is a factor of
another cyclic permutation $T_2$ if and only if the order of $T_1$ divides
the order of $T_2$.  In general, a cyclic permutation of order $n$ is a
factor of a transformation $T:X\to X$ if and only if the set $X$ can be
split into $n$ disjoint closed subsets $X_1,X_2,\dots,X_n$ which are
cyclically permuted by $T$, that is, $T(X_i)\subset X_{i+1}$ for
$1\le i\le n-1$ and $T(X_n)\subset X_1$.  Note that the sets
$X_1,X_2,\dots,X_n$ are both closed and open.

Now assume that a transformation $T:X\to X$ is a continuous extension of
two cyclic permutations $T_1$ and $T_2$ such that $T_1$ is a factor of
$T_2$.  Then there exist continuous onto maps $f:X_2\to X_1$,
$f_1:X\to X_1$, and $f_2:X\to X_2$ such that $fT_2=T_1f$, $f_1T=T_1f_1$,
and $f_2T=T_2f_2$, i.e., the following diagrams are commutative:
$$
\begin{array}{ccc}
X & \stackrel{T}{\longrightarrow} & X\\[0.2em]
f_2\Big\downarrow\phantom{f_2} && \phantom{f_2}\Big\downarrow f_2\\[0.25em]
X_2 & \stackrel{T_2}{\longrightarrow} & X_2\\[0.2em]
f\Big\downarrow\phantom{f} && \phantom{f}\Big\downarrow f\\[0.25em]
X_1 & \stackrel{T_1}{\longrightarrow} & X_1
\end{array}
\qquad\quad
\begin{array}{ccc}
X & \stackrel{T}{\longrightarrow} & X\\[0.2em]
f_1\Big\downarrow\phantom{f_1} && \phantom{f_1}\Big\downarrow f_1\\[0.25em]
X_1 & \stackrel{T_1}{\longrightarrow} & X_1
\end{array}
$$
Note that, given $f$ and $f_2$, we can always take $f_1=ff_2$.  It turns
out that the latter identity can also be satisfied when one is given $f$
and $f_1$ and has to choose $f_2$.

\begin{lemma}\label{odo1}
For any choice of the maps $f$ and $f_1$ above, we can choose the map
$f_2$ so that $ff_2=f_1$.
\end{lemma}

\begin{proof}
Take an arbitrary continuous map $h:X\to X_2$ such that $hT=T_2h$.  For any
$x\in X_1$ and $y\in X_2$ let $U(x)=f_1^{-1}(x)$, $V(x)=h^{-1}(y)$, and
$W(x,y)=U(x)\cap V(y)$.  Then $U(x)$, $x\in X_1$ is a collection of
disjoint closed sets that partition $X$.  The same holds true for the
collections $V(y)$, $y\in X_2$ and $W(x,y)$, $(x,y)\in X_1\times X_2$.
Furthermore, $T(U(x))\subset U(T_1x)$, $T(V(y))\subset V(T_2y)$, and
$T(W(x,y))\subset W(T_1x,T_2y)$.

By $n$ denote the cardinality of the set $X_1$.  For any $y\in X_2$ let
$Y(y)=\bigcup_{k=0}^{n-1}W(f(y),T_2^ky)$.  Clearly, each $Y(y)$ is a closed
subset of $U(f(y))$.  Besides, $T(Y(y))\subset Y(T_2y)$ since
$T_1f(y)=f(T_2y)$.  For any $x_0\in X_1$ and $y_0\in X_2$ there is a unique
$k\in\{0,1,\dots,n-1\}$ such that $f(T_2^{-k}y_0)=x_0$.  It follows that
the sets $Y(y)$, $y\in X_2$ are disjoint and cover the entire set $X$.

Define a map $f_2:X\to X_2$ so that $f_2(z)=y$ for all $y\in X_2$ and
$z\in Y(y)$.  Since $Y(y)$, $y\in X_2$ is a collection of disjoint closed
sets that partition $X$, the map $f_2$ is well defined and continuous.
Since $Y(y)\subset U(f(y))$ and $T(Y(y))\subset Y(T_2y)$ for any
$y\in X_2$, it follows that $ff_2=f_1$ and $f_2T=T_2f_2$.  Clearly, $f_2$
is onto.
\end{proof}

For any integer $n>0$ we denote $\bZ_n=\bZ/n\bZ$.  Given $k\in\bZ_n$ and
$m\in\bZ$, the sum $k+m$ is a well defined element of $\bZ_n$.  The {\em
odometer\/} on $\bZ_n$ is the transformation $x\mapsto x+1$.  Now let
$n_1,n_2,\dots$ be a finite or infinite sequence of positive integers.  The
{\em odometer\/} on $\bZ_{n_1}\times\bZ_{n_2}\times\dots$ is a
transformation $T$ defined as follows.  For any $m_i\in\bZ_{n_i}$,
$i=1,2,\dots$, we let $T(m_1,m_2,\dots)=(k_1,k_2,\dots)$, where $k_i=m_i+1$
if $m_j=-1+n_j\bZ$ for $1\le j<i$, and $k_i=m_i$ otherwise.  We regard each
$\bZ_n$ as a discrete topological space and endow $\bZ_{n_1}\times\bZ_{n_2}
\times\dots$ with the product topology.  Then the odometer is a
homeomorphism of a compact set.  It is easy to see that any odometer is
minimal.  The odometer on a finite set is a cyclic permutation.

\begin{lemma}\label{odo2}
For any odometer $T_0$ there exist a compact Abelian group $G$ and $g_0\in
G$ such that $T_0$ is continuously conjugated to the transformation
$g\mapsto g+g_0$ of $G$.
\end{lemma}

\begin{proof}
Suppose $T_0$ is the odometer on $X=\bZ_{n_1}\times\bZ_{n_2}\times\dots$.
If $X$ is finite then $T_0$ is a cyclic permutation, hence it is
conjugated to the odometer on some $\bZ_n$.  Now assume $X$ is
infinite.  Let us regard $\bZ$ as a discrete topological group and endow
the countable product $G=\bZ\times\bZ\times\dots$ with the product
topology.  By $G_0$ denote the closed subgroup of $G$ generated by elements
$(n_1,-1,0,0,0,\dots)$, $(0,n_2,-1,0,0,\dots)$, $(0,0,n_3,-1,0,\dots),
\dots$.  Let $g_0\in G/G_0$ be the coset containing $(1,0,0,\dots)$.  Each
$g\in G/G_0$ intersects the compact set $Y=\{0,1,\dots,n_1-1\}\times
\{0,1,\dots,n_2-1\}\times\dots$ in exactly one element.  For any
$y=(y_1,y_2,\dots)\in Y$ let $f_1(y)=(y_1+n_1\bZ,y_2+n_2\bZ,\dots)\in X$
and let $f_2(y)\in G/G_0$ be the coset containing $y$.  It is easy to
observe that the maps $f_1:Y\to X$ and $f_2:Y\to G/G_0$ are homeomorphisms.
Besides, $f_2f_1^{-1}(T_0x)=f_2f_1^{-1}(x)+g_0$ for all $x\in X$.
\end{proof}

\begin{lemma}\label{odo3}
Assume that a transformation $T:X\to X$ is a continuous extension of an
odometer $T_0:X_0\to X_0$.  Then for any $x\in X$ and $x_0\in X_0$ there
exists a continuous map $f:X\to X_0$ such that $f$ is onto, $fT=T_0f$, and
$f(x)=x_0$.  The map $f$ is unique provided that $T$ is transitive.
\end{lemma}

\begin{proof}
It is no loss of generality to replace $T_0$ by a continuously conjugated
transformation.  In view of Lemma \ref{odo2}, we can assume that $X_0$ is a
compact Abelian group and $T_0x_0=x_0+g_0$ for some $g_0\in X_0$ and all
$x_0\in X_0$.  Let $f:X\to X_0$ be a continuous onto map such that
$fT=T_0f$.  For any $g\in X_0$ and $x\in X$ let $f_g(x)=f(x)+g$.  Then
$f_g$ is a continuous map of $X$ onto $X_0$ and $f_gT(x)=fT(x)+g=T_0f(x)+g=
f(x)+g_0+g=T_0f_g(x)$ for all $x\in X$.  Obviously, for any $x\in X$ and
$x_0\in X_0$ there exists a unique $g\in X_0$ such that $f_g(x)=x_0$.

Now assume $T$ is transitive and pick $y\in X$ such that the orbit of $y$
under the transformation $T$ is dense in $X$.  Suppose $h:X\to X_0$ is a
continuous map such that $hT=T_0h$.  We have $h(y)=f_g(y)$ for some $g\in
X_0$.  Since $hT=T_0h$ and $f_gT=T_0f_g$, it follows that $h(T^iy)=
f_g(T^iy)$ for $i=1,2,\dots$.  Then density of the sequence $y,Ty,T^2y,
\dots$ in $X$ implies that $h=f_g$.
\end{proof}

\begin{lemma}\label{odo4}
Two odometers are continuously conjugated if either of them is a continuous
factor of the other.
\end{lemma}

\begin{proof}
Let $T_1:X_1\to X_1$ and $T_2:X_2\to X_2$ be odometers such that $T_1$ is
both a continuous factor and a continuous extension of $T_2$.  Pick
$x_1\in X_1$ and $x_2\in X_2$.  By Lemma \ref{odo3}, there are unique
continuous onto maps $f_1:X_1\to X_2$ and $f_2:X_2\to X_1$ such that
$f_1T_1=T_2f_1$, $f_2T_2=T_1f_2$, $f_1(x_1)=x_2$, and $f_2(x_2)=x_1$.  Note
that $f_1f_2f_1$ and $f_2f_1f_2$ are continuous maps, $f_1f_2f_1(X_1)=X_2$,
and $f_2f_1f_2(X_2)=X_1$.  Further, $(f_1f_2f_1)T_1=f_1f_2T_2f_1=
f_1T_1f_2f_1=T_2(f_1f_2f_1)$ and $f_1f_2f_1(x_1)=x_2$.  Similarly,
$(f_2f_1f_2)T_2=T_1(f_2f_1f_2)$ and $f_2f_1f_2(x_2)=x_1$.  It follows that
$f_1f_2f_1=f_1$ and $f_2f_1f_2=f_2$.  Since the maps $f_1$ and $f_2$ are
onto, $f_1f_2$ and $f_2f_1$ are the identity maps of $X_2$ and $X_1$,
respectively.  Thus $f_1$ and $f_2$ are homeomorphisms.
\end{proof}

The following two lemmas explore relations between odometers and cyclic
permutations.

\begin{lemma}\label{odo5}
Assume that a transformation $T$ is a continuous extension of cyclic
permutations of orders $n_1,n_1n_2,n_1n_2n_3,\dots$, where
$n_1,n_2,n_3,\dots$ is a sequence of positive integers.  Then $T$ is also
a continuous extension of the odometer on
$\bZ_{n_1}\times\bZ_{n_2}\times\bZ_{n_3}\times\dots$.
\end{lemma}

\begin{proof}
We assume that the sequence $n_1,n_2,n_3,\dots$ is infinite as otherwise
the lemma is trivial.  For any $k\ge1$ let $X_k=\bZ_{n_1}\times\bZ_{n_2}
\times\dots\times\bZ_{n_k}$ and denote by $T_k$ the odometer on $X_k$.
Also, let $T_\infty$ denote the odometer on $X_\infty=\bZ_{n_1}\times
\bZ_{n_2}\times\bZ_{n_3}\times\dots$ and let $X$ denote the space on which
$T$ acts.  For any $k\ge1$ consider the natural projections
$\pi_k:X_{k+1}\to X_k$ and $p_k:X_\infty\to X_k$.  They are continuous and
onto.  Besides, $\pi_kT_{k+1}=T_k\pi_k$ and $p_kT_\infty=T_kp_k$.
Since the odometer $T_k$ is a cyclic permutation of order
$n_1n_2\cdots n_k$, it is a factor of the transformation $T$.  Hence there
is a continuous map $f_k:X\to X_k$ such that $f_kT=T_kf_k$.  In view of
Lemma \ref{odo1}, we can choose the maps $f_1,f_2,\dots$ so that
$f_k=\pi_kf_{k+1}$ for all $k\ge1$.

Define a map $f:X\to X_\infty$ as follows.  Given $x\in X$, let
$f(x)=(m_1,m_2,\dots)$, where $f_k(x)=(m_1,m_2,\dots,m_k)$ for all $k\ge1$.
The map $f$ is well defined since $f_k=\pi_kf_{k+1}$ for all $k\ge1$.  Its
continuity follows from the continuity of $f_1,f_2,\dots$.  Furthermore,
$fT=T_\infty f$ as $f_kT=T_kf_k$ and $p_kT_\infty=T_kp_k$ for all $k\ge1$.
In particular, the image $f(X)$ is invariant under $T_\infty$.  Since $X$
is compact, $f(X)$ a nonempty compact subset of $X_\infty$.  The minimality
of the odometer $T_\infty$ implies that the map $f$ is onto.  Thus
$T_\infty$ is a continuous factor of $T$.
\end{proof}

\begin{lemma}\label{odo6}
A cyclic permutation of order $n$ is a continuous factor of the odometer on
$\bZ_{m_1}\times\bZ_{m_2}\times\bZ_{m_3}\times\dots$ if and only if $n$
divides some of the numbers $m_1,m_1m_2,m_1m_2m_3,\dots$.
\end{lemma}

\begin{proof}
Let $T$ be the odometer on $X=\bZ_{m_1}\times\bZ_{m_2}\times\dots$.  First
suppose the sequence $m_1,m_2,\dots$ is finite.  Denote by $k$ its length.
Then $T$ is a cyclic permutation of order $m_1m_2\cdots m_k$.  Hence a
cyclic permutation of order $n$ is a factor of $T$ if and only if $n$
divides $m_1m_2\cdots m_k$.

Now consider the case when the sequence $m_1,m_2,\dots$ is infinite.
Suppose that a cyclic permutation $T_0:X_0\to X_0$ of order $n$ is a factor
of $T$.  Let $f:X\to X_0$ be a continuous map such that $fT=T_0f$.  It is
easy to observe that for any $x\in X$ the sequence
$T^{m_1}x,T^{m_1m_2}x,T^{m_1m_2m_3}x,\dots$ converges to $x$.  Hence the
sequence $f(T^{m_1}x),f(T^{m_1m_2}x),f(T^{m_1m_2m_3}x),\dots$ converges to
$f(x)$.  Since $X_0$ is a finite set with the discrete topology, this means
that $f(T^{m_1m_2\cdots m_k}x)=T_0^{m_1m_2\cdots m_k}f(x)$ coincides with
$f(x)$ for large $k$.  It follows that $n$ divides $m_1m_2\cdots m_k$ for
large $k$.

Conversely, if $n$ divides some $m_1m_2\cdots m_k$ then any cyclic
permutation $T_0$ of order $n$ is a factor of the odometer $T_k$ on
$X_k=\bZ_{m_1}\times\bZ_{m_2}\times\dots\times\bZ_{m_k}$.  But $T_k$ is
a factor of $T$ since the natural projection $p_k:X\to X_k$ is continuous
and satisfies $p_kT=T_kp_k$.  Then $T_0$ is also a factor of the odometer
$T$.
\end{proof}

\begin{lemma}\label{odo7}
The odometer on $\bZ_{n_1}\times\bZ_{n_2}\times\bZ_{n_3}\times\dots$ is a
continuous factor of the odometer on $\bZ_{m_1}\times\bZ_{m_2}\times
\bZ_{m_3}\times\dots$ if and only if each element of the sequence
$n_1,n_1n_2,n_1n_2n_3,\dots$ divides an element of the sequence
$m_1,m_1m_2,m_1m_2m_3,\dots$.
\end{lemma}

\begin{proof}
Let $T_1$ denote the odometer on $\bZ_{n_1}\times\bZ_{n_2}\times\bZ_{n_3}
\times\dots$ and $T_2$ denote the odometer on $\bZ_{m_1}\times\bZ_{m_2}
\times\bZ_{m_3}\times\dots$.  By Lemma \ref{odo6}, cyclic permutations of
orders $n_1,n_1n_2,n_1n_2n_3,\dots$ are factors of $T_1$.  Assume that
$T_1$ is a factor of $T_2$.  Then all factors of $T_1$ are also factors of
$T_2$.  It follows from Lemma \ref{odo6} that each of the numbers
$n_1,n_1n_2,n_1n_2n_3,\dots$ divides some of the numbers
$m_1,m_1m_2,m_1m_2m_3,\dots$.

Conversely, assume that each of the numbers $n_1,n_1n_2,n_1n_2n_3,\dots$
divides some of the numbers $m_1,m_1m_2,m_1m_2m_3,\dots$.  Then Lemma
\ref{odo6} implies that cyclic permutations of orders
$n_1,n_1n_2,n_1n_2n_3,\dots$ are factors of the odometer $T_2$.  By Lemma
\ref{odo5}, the odometer $T_1$ is a factor of $T_2$ as well.
\end{proof}

To any continuous transformation $T$ of a compact topological space we
associate the set $\CF(T)$ of positive integers $n$ such that $T$ is a
continuous extension of the cyclic permutation of order $n$ ($\CF$ stands
for ``cyclic factors'').

\begin{lemma}\label{odo8}
The set $\CF(T)$ has the following properties:

(i) $1\in \CF(T)$;

(ii) if $n\in\CF(T)$ and $d>0$ is a divisor of $n$, then $d\in\CF(T)$;

(iii) if $n_1,n_2,\dots,n_k\in\CF(T)$, then
$\lcm(n_1,n_2,\dots,n_k)\in\CF(T)$.
\end{lemma}

\begin{proof}
Property (i) is trivial.

Suppose $d$ and $n$ are positive integers.  If $d$ divides $n$ then
a cyclic permutation of order $d$ is a factor of a cyclic permutation of
order $n$.  Therefore all continuous extensions of the latter permutation
are also continuous extensions of the former one.  In particular,
$d\in\CF(T)$ whenever $n\in\CF(T)$.  Property (ii) is verified.

Let $X$ denote the topological space on which $T$ acts.  Given
$m,n\in\CF(T)$, there exist continuous maps $f_1:X\to\bZ_m$ and
$f_2:X\to\bZ_n$ such that $f_1(Tx)=f_1(x)+1$ and $f_2(Tx)=f_2(x)+1$ for all
$x\in X$.  Then $f=(f_1,f_2)$ is a continuous map of $X$ to
$\bZ_m\times\bZ_n$.  Furthermore, $fT=T_0f$, where $T_0$ denotes the
transformation $(x,y)\mapsto(x+1,y+1)$ of $\bZ_m\times\bZ_n$.  Assume that
$m$ and $n$ are coprime.  Then $T_0$ is a cyclic permutation of order $mn$
and the map $f$ is onto.  Consequently, $mn\in\CF(T)$.

Let $n_1,n_2\in\CF(T)$.  It is easy to show that there exist positive
coprime integers $d_1,d_2$ such that $d_1$ divides $n_1$, $d_2$ divides
$n_2$, and $d_1d_2=\lcm(n_1,n_2)$.  By the above $d_1,d_2,d_1d_2\in\CF(T)$.
Now property (iii) follows by induction.
\end{proof}

The continuous cyclic factors of a transformation $T:X\to X$ are related to
continuous eigenfunctions of $T$.  Let $\cU_T$ denote a linear operator
that acts on functions on $X$ by precomposing them with $T$: $\cU_T\phi=
\phi T$.  A nonzero function $\phi:X\to\bC$ is an {\em eigenfunction\/} of
$T$ associated with an {\em eigenvalue\/} $\la$ if $\cU_T\phi=\la\phi$,
that is, if $\phi(Tx)=\la\phi(x)$ for all $x\in X$.

\begin{lemma}\label{eigen1}
If $n\in\CF(T)$, then the transformation $T$ admits a continuous
eigenfunction associated with the eigenvalue $e^{2\pi i/n}$, a primitive
$n$th root of unity.  For a transitive $T$, the converse is true as well.
Besides, for a transitive $T$ any continuous eigenfuction is determined by
its eigenvalue uniquely up to scaling.
\end{lemma}

\begin{proof}
Given a positive integer $n$, let $R_n$ denote the set of all $n$th roots
of unity in $\bC$.  $R_n$ is a multiplicative cyclic group of order $n$
generated by a primitive root $\zeta_n=e^{2\pi i/n}$.  A transformation
$T_n:R_n\to R_n$ defined by $T_n(z)=\zeta_nz$ is a cyclic permutation of
order $n$.  If $n\in\CF(T)$ for some transformation $T:X\to X$, then there
exists a continuous mapping $f:X\to R_n$ such that $fT=T_nf$, i.e.,
$f(Tx)=T_n(f(x))=\zeta_nf(x)$ for all $x\in X$.  Clearly, $f$ is a
continuous eigenfunction of $T$ associated with the eigenvalue $\zeta_n$.

Now assume that $T:X\to X$ is transitive and pick a point $x_0\in X$ with
dense orbit.  Let $\phi$ be a continuous eigenfunction of $T$ associated
with an eigenvalue $\la$.  The function $\phi$ is uniquely determined by
its values on the orbit $x_0,Tx_0,T^2x_0,\dots$.  Since $\phi(T^kx_0)=
\la^k\phi(x_0)$ for $k=1,2,\dots$, the function $\phi$ is uniquely
determined by $\phi(x_0)$ and $\la$.  Note that $\phi(x_0)\ne0$ as
otherwise $\phi$ will be identically zero.  Since nonzero scalar multiples
of $\phi$ are also eigenfunctions with the same eigenvalue, $\phi$ is
determined by the eigenvalue $\la$ up to scaling.

In the case $\la=\zeta_n$, replace the eigenfunction $\phi$ by a scalar
multiple so that $\phi(x_0)=1$.  Then $\phi(T^kx_0)\in R_n$ for all
$k\ge1$, which implies that $\phi(x)\in R_n$ for all $x\in X$.  Therefore
$\phi$ maps $X$ onto $R_n$ and $\phi(Tx)=\zeta_n\phi(x)=T_n(\phi(x))$ for
all $x\in X$.  Thus the cyclic permutation $T_n$ is a factor of $T$, that
is, $n\in\CF(T)$.
\end{proof}

\begin{lemma}\label{odo9}
Suppose $F_0$ is a set of positive integers such that (i) $1\in F_0$,
(ii) any positive divisor of an element of $F_0$ also belongs to $F_0$, and
(iii) the least common multiple of finitely many elements of $F_0$ is in
$F_0$ as well.
Then there exists an odometer $T_0$ such that $\CF(T_0)=F_0$.
\end{lemma}

\begin{proof}
First suppose that the set $F_0$ is finite.  Let $N$ be its maximal
element.  For any $m\in F_0$ the least common multiple of $m$ and $N$
belongs to $F_0$.  By the choice of $N$, we have $\lcm(m,N)=N$, that is,
$m$ divides $N$.  On the other hand, $F_0$ contains all positive divisors
of $N$.  Therefore $F_0$ is the set of positive integers that divide $N$.
It follows that $F_0=\CF(T_0)$ for any cyclic permutation $T_0$ of order
$N$.  One example of such a permutation is the odometer on $\bZ_N$.

Now suppose that the set $F_0$ is infinite and let $n_1,n_2,\dots$ be a
complete list of its elements.  For any $k\ge1$ let $m_k$ be the least
common multiple of the integers $n_1,n_2,\dots,n_k$.  Then the numbers
$m_1,m_2,\dots$ belong to $F_0$, each $m_k$ divides $m_{k+1}$, and any
$n\in F_0$ divides some $m_k$.  On the other hand, all positive divisors of
any $m_k$ are in $F_0$.  Hence $F_0$ is the set of positive integers that
divide some of the numbers $m_1,m_2,\dots$.  It follows from Lemma
\ref{odo6} that $F_0=\CF(T_0)$, where $T_0$ is the odometer on
$\bZ_{m_1}\times\bZ_{m_2/m_1}\times\bZ_{m_3/m_2}\times\dots$.
\end{proof}

\begin{proposition}\label{odo10}
Suppose $T_1$ and $T_2$ are odometers. Then

(i) $T_1$ is a continuous factor of $T_2$ if and only if $\CF(T_1)\subset
\CF(T_2)$;

(ii) $T_1$ and $T_2$ are continuously conjugated if and only if
$\CF(T_1)=\CF(T_2)$.
\end{proposition}

\begin{proof}
Let $T_1$ be the odometer on $\bZ_{n_1}\times\bZ_{n_2}\times\bZ_{n_3}\times
\dots$ and $T_2$ be the odometer on $\bZ_{m_1}\times\bZ_{m_2}\times
\bZ_{m_3}\times\dots$.  Lemma \ref{odo6} implies that $\CF(T_1)$ is the set
of positive integers that divide some of the numbers
$n_1,n_1n_2,n_1n_2n_3,\dots$.  Likewise, $\CF(T_2)$ is the set of positive
integers that divide some of the numbers $m_1,m_1m_2,m_1m_2m_3,\dots$.
Therefore $\CF(T_1)\subset\CF(T_2)$ if and only if each of the numbers
$n_1,n_1n_2,n_1n_2n_3,\dots$ divides some of $m_1,m_1m_2,m_1m_2m_3,\dots$.
According to Lemma \ref{odo7}, this is exactly when $T_1$ is a factor of
$T_2$.

By Lemma \ref{odo4}, the odometers $T_1$ and $T_2$ are conjugated if either
of them is a factor of the other.  Hence the statement (ii) of the
proposition follows from the statement (i).
\end{proof}

Suppose that an odometer $T_0$ is a continuous factor of a transformation
$T$.  We shall say that $T_0$ is a {\em maximal odometer factor\/} of $T$
if any odometer is a continuous factor of $T_0$ whenever it is a continuous
factor of $T$.

\begin{proposition}\label{odo11}
(i) An odometer $T_0$ is a continuous factor of a transformation $T$
if and only if $\CF(T_0)\subset\CF(T)$.

(ii) An odometer $T_0$ is a maximal odometer factor of $T$ if and only
if $\CF(T_0)=\CF(T)$.

(iii) The maximal odometer factor always exists and is unique up to
continuous conjugacy.
\end{proposition}

\begin{proof}
Let $T_0$ be the odometer on $\bZ_{n_1}\times\bZ_{n_2}\times\bZ_{n_3}\times
\dots$.  By Lemma \ref{odo6}, the integers $n_1,n_1n_2,n_1n_2n_3,\dots$
belong to $\CF(T_0)$.  If $\CF(T_0)\subset\CF(T)$, they belong to $\CF(T)$
as well.  Then it follows from Lemma \ref{odo5} that $T_0$ is a factor of
$T$.  Conversely, if $T_0$ is a factor of $T$ then all factors of $T_0$ are
also factors of $T$; in particular, $\CF(T_0)\subset\CF(T)$.

Lemmas \ref{odo8} and \ref{odo9} imply that for any transformation $T$
there exists an odometer $T_1$ such that $\CF(T_1)=\CF(T)$.  By the above
$T_1$ is a factor of $T$.  Suppose $T_2$ is another odometer that is a
factor of $T$.  Then $\CF(T_2)\subset\CF(T)=\CF(T_1)$.  By Proposition
\ref{odo10}, $T_2$ is a factor of $T_1$.  Therefore $T_1$ is a maximal
odometer factor of $T$.  Its uniqueness up to continuous conjugacy follows
from Lemma \ref{odo4}.  Thus an odometer $T_0$ is a maximal odometer factor
of $T$ if and only if it is conjugated to $T_1$.  According to Proposition
\ref{odo10}, this is exactly when $\CF(T_0)=\CF(T_1)=\CF(T)$.
\end{proof}

Now let us consider ergodic properties of odometers.  First recall some
definitions.  Let $X$ be a measured space and $T:X\to X$ be a measurable
transformation.  A measure $\mu$ on $X$ is {\em invariant\/} under $T$ if
$\mu(T^{-1}(B))=\mu(B)$ for any measurable subset $B\subset X$.  Let
$\cU_T$ be the linear operator acting on functions on $X$ by precomposition
with $T$: $\cU_T\phi=\phi T$.  If $\mu$ is an invariant measure, then
$\cU_T$ is a unitary operator when restricted to the Hilbert space
$L_2(X,\mu)$.  The spectral properties of this unitary operator are
referred to as {\em spectral properties\/} of the dynamical system
$(X,\mu,T)$.  For instance, one says that $T$ has {\em pure point
spectrum\/} if $L_2(X,\mu)$ admits an orthonormal basis consisting of
eigenfunctions of $\cU_T$.

The transformation $T$ is {\em ergodic\/} with respect to the invariant
measure $\mu$ if for any measurable subset $B\subset X$ that is backward
invariant under $T$ (i.e., $T^{-1}(B)\subset B$) one has $\mu(B)=0$ or
$\mu(X\setminus B)=0$.  If $T$ is ergodic and the measure $\mu$ is finite,
then Birkhoff's ergodic theorem implies that $\mu$-almost all orbits of $T$
are uniformly distributed in $X$ relative to this measure.

A homeomorphism $T$ of a compact topological space $X$ is called
{\em uniquely ergodic\/} if there exists a unique Borel probability measure
on $X$ invariant under $T$.  The uniquely ergodic transformation $T$ is
ergodic with respect to the unique invariant measure.  Moreover, in this
case every orbit of $T$ is uniformly distributed in $X$.

\begin{proposition}\label{eigen2}
Let $T_0$ be an odometer.  Then

(i) $T_0$ is uniquely ergodic and has purely point spectrum;

(ii) the eigenvalues of $T_0$ are all $n$th roots of unity, where $n$ runs
through the set $\CF(T_0)$;

(iii) all eigenvalues of $T_0$ are simple and all eigenfunctions are
continuous.
\end{proposition}

\begin{proof}
Let $T_0$ be the odometer on $X=\bZ_{n_1}\times\bZ_{n_2}\times\dots$.
Denote by $\mu_0$ a Borel probability measure on $X$ that is the direct
product of normalized counting measures on finite sets $\bZ_{n_1},\bZ_{n_2},
\dots$.  For any cylindrical set $C$ of the form $\{z\}\times\bZ_{n_{k+1}}
\times\bZ_{n_{k+1}}\times\dots$, where $z\in\bZ_{n_1}\times\dots\times
\bZ_{n_k}$, we have $\mu_0(C)=(n_1n_2\cdots n_k)^{-1}$.  For a fixed $k$,
there are $n_1n_2\cdots n_k$ such sets.  They partition the set $X$ and are
cyclically permuted by the odometer $T_0$.  This implies that
$\mu_0(T_0^{-1}(C))=\mu_0(C)$ for all cylindrical sets $C$.  Also,
$\mu(C)=\mu_0(C)$ for any Borel probability measure on $X$ invariant under
$T_0$.  Since any Borel measure on $X$ is determined by its values on
cylindrical sets, we obtain that the odometer $T_0$ is uniquely ergodic and
$\mu_0$ is the unique invariant measure.

Let $E$ be the set of all $n$th roots of unity, where $n$ runs through the
set $\CF(T_0)$.  Take any $\zeta\in E$.  We have $\zeta^n=1$ for some
$n\in\CF(T_0)$.  Then $\zeta=\zeta_n^k$, where $\zeta_n=e^{2\pi i/n}$ and
$k$ is a positive integer.  By Lemma \ref{eigen1}, the odometer $T_0$
admits a continuous eigenfunction $f_n$ associated with the eigenvalue
$\zeta_n$.  It is easy to observe that $f_n^k$ is a continuous
eigenfunction of $T_0$ associated with the eigenvalue $\zeta$.  According
to Lemma \ref{eigen1}, any continuous eigenfuction of $T_0$ is determined
by its eigenvalue uniquely up to scaling.

To finish the proof of the proposition, it remains to show that the set $F$
of continuous eigenfunctions of $T_0$ associated with eigenvalues from the
set $E$ is complete in the Hilbert space $L_2(X,\mu_0)$, i.e., the linear
span of $F$ is dense in $L_2(X,\mu_0)$.  Take an arbitrary cylindrical set
$C=\{z\}\times\bZ_{n_{k+1}}\times\bZ_{n_{k+1}}\times\dots$.  The product
$m=n_1n_2\cdots n_k$ belongs to $\CF(T_0)$ due to Lemma \ref{odo6}.  Then a
primitive $m$th root of unity $\zeta_m=e^{2\pi i/m}$ belongs to $E$.  Let
$f$ be a continuous eigenfunction of $T_0$ associated with the eigenvalue
$\zeta_m$.  The minimality of $T_0$ implies that $f$ is nowhere zero.
Replacing $f$ by a scalar multiple, we can assume that $f(x_0)=1$ for some
$x_0\in C$.  Then $f(T_0^kx_0)=\zeta_m^k$ for $k=1,2,\dots$.  Note that
$T_0^kx_0\in C$ if and only if $k$ divides $m$.  Besides, $f(T_0^kx_0)=1$
if and only if $k$ divides $m$.  It follows that all values of $f$ are
$m$th roots of unity, moreover, $f(x)=1$ if and only if $x\in C$.
Therefore the sum $f+f^2+\cdots+f^m$ is identically $m$ on the set $C$ and
identically zero elsewhere.  Notice that each term in this sum is a
continuous eigenfunction of $T_0$ with eigenvalue an $m$th root of unity.
Thus the characteristic functions of cylindrical sets are contained in the
span of $F$.  These characteristic functions form a complete set in
$L_2(X,\mu)$ for any finite Borel measure $\mu$ on $X$.
\end{proof}

Sometimes the odometers as defined above in this section are called {\em
generalized odometers\/} while the notion ``odometer'' refers to $p$-adic
odometers, which are defined as follows.
Let $p$ be a prime integer and $0<\rho<1$.  Any nonzero $r\in\bQ$ is
uniquely represented in the form $p^k\frac{m}{n}$, where $k,m,n$ are
integers, $n>0$, and $p$ divides neither $m$ nor $n$.  We let $|r|_p=
\rho^k$.  Also, let $|0|_p=0$.  Now $|\cdot|_p$ is a norm on the field
$\bQ$ called the {\em $p$-adic norm}.  The $p$-adic norm induces a distance
$d_p$ on $\bQ$, $d_p(r_1,r_2)=|r_1-r_2|_p$ for all $r_1,r_2\in\bQ$.  By
definition, the field $\cF_p$ of {\em $p$-adic numbers\/} is the completion
of the field $\bQ$ with respect to the $p$-adic norm.  The ring $\cZ_p$ of
{\em $p$-adic integers\/} is the closure of the ring $\bZ$ in $\cF_p$.  The
transformation $x\mapsto x+1$ of $\cZ_p$ is called the {\em $p$-adic
odometer}.  The $2$-adic odometer is also called the {\em binary odometer}.

\begin{lemma}\label{odo12}
The $p$-adic odometer is continuously conjugated to the odometer on
$\bZ_p\times\bZ_p\times\bZ_p\times\dots$.
\end{lemma}

\begin{proof}
Let $X$ denote $\bZ_p\times\bZ_p\times\dots$ and $T$ denote the odometer on
$X$.  By $X_0$ denote the countable product
$\{0,1,\dots,p-1\}\times\{0,1,\dots,p-1\}\times\dots$.  We endow $X_0$ with
the product topology.  Then the map $f:X_0\to X$ defined by
$f(n_1,n_2,\dots)=(n_1+p\bZ,n_2+p\bZ,\dots)$ is a homeomorphism.

An arbitrary $p$-adic integer is uniquely expanded into a series of the
form $\sum_{k=1}^\infty n_kp^{k-1}$, where $n_k\in\{0,1,\dots,p-1\}$ for
$k=1,2,\dots$.  Moreover, any series of this form converges in $\cZ_p$.  It
follows that the map $h:X_0\to\cZ_p$ defined by $h(n_1,n_2,\dots)=
\sum_{k=1}^\infty n_kp^{k-1}$ is one-to-one and onto.  It is easy to see
that $h$ is continuous as well.  Since $X_0$ is compact, $h$ is a
homeomorphism.

It follows from the definition of the maps $f$ and $h$ that
$hf^{-1}(Tx)=hf^{-1}(x)+1$ for all $x\in X$.  Thus $T$ is continuously
conjugated to the $p$-adic odometer.
\end{proof}

\section{Toeplitz sequences}\label{toe}

Let $\cA$ be a nonempty finite set.  We denote by $\cA^{\bN}$ the countable
product $\cA\times\cA\times\cA\times\dots$ endowed with the product
topology (here $\bN$ refers to positive integers).  Any $\om\in\cA^{\bN}$
is represented as an infinite sequence $(\om_1,\om_2,\dots)$ of elements of
$\cA$.  Denote by $\si$ the map on $\cA^{\bN}$ that sends any sequence
$\om$ to the sequence obtained by deleting the first element of $\om$.
That is, if $\om=(\om_1,\om_2,\dots)$ then $\si\om=(\om'_1,\om'_2,\dots)$,
where $\om'_k=\om_{k+1}$ for $k=1,2,\dots$.  The map $\si$ is called the
{\em (one-sided) shift}.  It is a continuous map of the compact topological
space $\cA^{\bN}$ onto itself.  The restriction of the shift to any
closed invariant subset of $\cA^{\bN}$ is called a {\em subshift}.
In particular, any sequence $\om\in\cA^{\bN}$ gives rise to the subshift
$\si|_{Z(\om,\si)}$, where $Z(\om,\si)$ is the closure of the orbit
$\om,\si\om,\si^2\om,\dots$.

A sequence $\om=(\om_1,\om_2,\dots)\in\cA^{\bN}$ is called a {\em Toeplitz
sequence\/} if for any positive integer $n$ there exists a positive integer
$p$ such that $\om_n=\om_{n+kp}$ for $k=1,2,\dots$.  If $\om$ is a Toeplitz
sequence, then all shifted sequences $\si\om,\si^2\om,\dots$ are also
Toeplitz sequences.

\begin{lemma}\label{toe1}
Let $\om=(\om_1,\om_2,\dots)$ be a Toeplitz sequence.  Suppose $n$ and $p$
are positive integers such that $\om_n=\om_{n+kp}$ for $k=1,2,\dots$.
If $p<n$ then $\om_{n-p}=\om_n$.
\end{lemma}

\begin{proof}
Since $\om$ is a Toeplitz sequence, there exists a positive integer $q$
such that $\om_{n-p}=\om_{n-p+kq}$ for $k=1,2,\dots$.  In particular,
$\om_{n-p}=\om_{n-p+pq}$.  On the other hand, $\om_{n-p+pq}=\om_{n+(q-1)p}=
\om_n$.
\end{proof}

\begin{lemma}\label{toe2}
For any Toeplitz sequence $\om$ the subshift $\si|_{Z(\om,\si)}$ is
minimal.
\end{lemma}

\begin{proof}
For any integer $n\ge1$ there exists an integer $p_n\ge1$ such that
$\om_n=\om_{n+kp_n}$ for $k=1,2,\dots$.  Take an arbitrary integer $N\ge1$
and let $p=p_1p_2\dots p_N$.  Then $\om_n=\om_{n+kp}$ for $1\le n\le N$ and
$k\ge1$.  In particular, for any $k\ge1$ the first $N$ elements of the
sequence $\si^{kp}\om$ are the same as the first $N$ elements of $\om$.

Given $\om'\in Z(\om,\si)$, there are nonnegative integers $n_1,n_2,\dots$
such that $\si^{n_k}\om\to\om'$ as $k\to\infty$.  It is no loss to assume
that the numbers $n_1,n_2,\dots$ have the same remainder $r$ under division
by $p$.  Then each $\si^{n_k+p-r}\om$ has the same first $N$ elements as
$\om$.  Since $\si^{n_k+p-r}\om\to\si^{p-r}\om'$ as $k\to\infty$, the
sequence $\si^{p-r}\om'$ also has the same first $N$ elements as $\om$.
Since $N$ can be chosen arbitrarily large, it follows that $\om$ is in the
closure of the orbit $\om',\si\om',\si^2\om',\dots$.  Hence the orbit
$\om',\si\om',\si^2\om',\dots$ is dense in $Z(\om,\si)$.
\end{proof}

\begin{lemma}\label{toe3}
Suppose $\om$ is a Toeplitz sequence.  Then for any integer $p>0$ there
exists an integer $K>0$ with the following property.  If for some
$\om'=(\om'_1,\om'_2,\dots)\in Z(\om,\si)$ and integer $n>0$ we have that
$\om'_n=\om'_{n+kp}$ for $1\le k\le K$, then $\om'_n=\om'_{n+kp}$ for all
$k\ge1$.
\end{lemma}

\begin{proof}
For any $K\ge1$ let $S(p,K)$ denote the set of positive integers $n$ such
that $\om_n\ne\om_{n+kp}$ for some $1\le k\le K$.  By $P(p)$ denote the set
of positive integers $n$ such that $\om_n=\om_{n+kp}$ for all $k\ge1$.
Take any $q\in\{1,2,\dots,p\}$.  If $q\in P(p)$ then the numbers
$q+p,q+2p,\dots$ are in $P(p)$ as well.  Now suppose that $q\notin P(p)$.
Then $\om_{q+kp}\ne\om_q$ for some $k\ge1$.  Since $\om$ is a Toeplitz
sequence, there exist $p_1,p_2>0$ such that $\om_{q+mp_1}=\om_q$ and
$\om_{q+kp+mp_2}=\om_{q+kp}$ for all $m\ge1$.  In particular,
$\om_{q+lp_1p_2p}\ne\om_{q+kp+lp_1p_2p}$ for $l=0,1,2,\dots$.  It follows
that each of the numbers $q,q+p,q+2p,\dots$ belongs to the set
$S(p,k+p_1p_2)$.

Let $q\in\{1,2,\dots,p\}$.  By the above the numbers $q,q+p,q+2p,\dots$
either all belong to $P(p)$, or else they all belong to $S(p,K)$ for some
$K\ge1$.  Hence there exists an integer $K_0\ge1$ such that
$\bN=P(p)\cup S(p,K_0)$.

Let $\om'=(\om'_1,\om'_2,\dots)\in Z(\om,\si)$.  Suppose that for some
integer $n\ge1$ we have that $\om'_n=\om'_{n+kp}$ for $1\le k\le K_0$.
Since $\om'\in Z(\om,\si)$, there are nonnegative integers $n_1,n_2,\dots$
such that $\si^{n_m}\om\to\om'$ as $m\to\infty$.  If $m$ is large enough,
then the first $n+K_0p$ elements of the sequence $\si^{n_m}\om$ are the
same as the first $n+K_0p$ elements of $\om'$.  In particular,
$\om_{n_m+n}=\om_{n_m+n+kp}$ for $1\le k\le K_0$.  This means that
$n_m+n\notin S(p,K_0)$.  As $\bN=P(p)\cup S(p,K_0)$, we obtain that
$n_m+n\in P(p)$, i.e., $\om_{n_m+n}=\om_{n_m+n+kp}$ for all $k\ge1$.
Since $\si^{n_m}\om\to\om'$ as $m\to\infty$, it follows that
$\om'_n=\om'_{n+kp}$ for all $k\ge1$.
\end{proof}

\begin{lemma}\label{toe4}
Suppose $\om$ is a Toeplitz sequence.  Then the Toeplitz sequences in
$Z(\om,\si)$ form a residual (dense $G_\delta$) subset.
\end{lemma}

\begin{proof}
Given any positive integers $n,p,k$, let $T(n,p,k)$ be the set of all
sequences $\om'=(\om'_1,\om'_2,\dots)\in\cA^{\bN}$ such that
$\om'_n=\om'_{n+ip}$ for $i=1,2,\dots,k$.  Clearly, $T(n,p,k)$ is an open
subset of $\cA^{\bN}$.  The set $\cT$ of all Toeplitz sequences in
$\cA^{\bN}$ can be represented as
$$
\cT=\bigcap_{n=1}^\infty\bigcup_{p=1}^\infty\bigcap_{k=1}^\infty T(n,p,k).
$$

Suppose $\om$ is a Toeplitz sequence.  According to Lemma \ref{toe3}, for
any $p\ge1$ there exists an integer $K_p\ge1$ such that
$$
\bigcap_{k=1}^\infty T(n,p,k)\cap Z(\om,\si)=T(n,p,K_p)\cap Z(\om,\si)
$$
for all $n\ge1$.  Then
$$
\cT\cap Z(\om,\si)=\bigcap_{n=1}^\infty
\biggl(\, \bigcup_{p=1}^\infty T(n,p,K_p) \,\biggr) \cap Z(\om,\si).
$$
Since $\bigcup_{p=1}^\infty T(n,p,K_p)$ is an open subset of $\cA^{\bN}$,
it follows that $\cT\cap Z(\om,\si)$ is a $G_\delta$ subset of
$Z(\om,\si)$.  It is dense in $Z(\om,\si)$ since
$\om,\si\om,\si^2\om,\dots$ are Toeplitz sequences.
\end{proof}

Let $\om=(\om_1,\om_2,\dots)$ be a Toeplitz sequence.  We shall say that a
positive integer $p$ is a {\em partial period\/} of $\om$ if there is
$n\ge1$ such that $\om_n=\om_{n+kp}$ for $k=1,2,\dots$.
The integer $p$ is called an {\em essential partial period\/} of $\om$ if
there exists $m\ge1$ such that $\om_m=\om_{m+kp}$ for $k=1,2,\dots$ while
for any $1\le p'<p$ the sequence $\om_m,\om_{m+p'},\om_{m+2p'},\dots$
contains an element different from $\om_m$.
We denote by $\EP(\om)$ the set of all essential partial periods of $\om$
($\EP$ stands for ``essential periods'').

Recall from Section \ref{odo} that to each continuous self-mapping $T$ of
a compact topological space we associate the set $\CF(T)$ of positive
integers $n$ such that the cyclic permutation of order $n$ is a factor of
$T$.  It turns out that essential partial periods of a Toeplitz sequence
$\om$ completely determine the set $\CF(\si|_{Z(\om,\si)})$.

\begin{lemma}\label{toe5}
For any Toeplitz sequence $\om$ one has $\EP(\om)\subset
\CF(\si|_{Z(\om,\si)})$.
\end{lemma}

\begin{proof}
Let $\om=(\om_1,\om_2,\dots)$ be a Toeplitz sequence and $p\in\EP(\om)$.
Pick $m\ge1$ such that $\om_m=\om_{m+kp}$ for $k=1,2,\dots$ while
for any $1\le p'<p$ the sequence $\om_m,\om_{m+p'},\om_{m+2p'},\dots$
contains an element different from $\om_m$.

Denote $\om_m$ by $a$.  For any $\om'=(\om_1,\om_2,\dots)\in Z(\om,\si)$
consider the set $R(\om')$ of cosets $\alpha\in\bZ/p\bZ$ such that
$\om'_n=a$ for all $n\in\alpha$, $n>0$.  Lemma \ref{toe1} implies that
$m+p\bZ\in R(\om)$.  Besides, $p'+R(\om)\ne R(\om)$ for any $1\le p'<p$ as
otherwise $\om_m=\om_{m+p'}=\om_{m+2p'}=\dots$  It follows that the sets
$R(\om),1+R(\om),\dots,(p-1)+R(\om)$ are all distinct.

By Lemma \ref{toe1}, $R(\si\om')=-1+R(\om')$ for all $\om'\in Z(\om,\si)$.
By Lemma \ref{toe3}, the set $R(\om')$ is locally constant as a function of
$\om'$.  It follows that for any $\om'\in Z(\om,\si)$ we have
$R(\om')=-\alpha(\om')+R(\om)$, where $\alpha(\om')\in\bZ/p\bZ$.  Moreover,
$\alpha(\om')$ is uniquely determined by $\om'$, the mapping
$\om'\mapsto\alpha(\om')$ is continuous, and
$\alpha(\si\om')=\alpha(\om')+1$.  Thus the odometer on $\bZ/p\bZ$ is a
continuous factor of the subshift $\si|_{Z(\om,\si)}$.  That is,
$p\in\CF(\si|_{Z(\om,\si)})$.
\end{proof}

\begin{proposition}\label{toe6}
Suppose $\om$ is a Toeplitz sequence.  Then a positive integer $n$ belongs
to $\CF(\si|_{Z(\om,\si)})$ if and only if $n$ divides the least common
multiple of some $n_1,n_2\dots,n_k\in\EP(\om)$.
Equivalently, $\CF(\si|_{Z(\om,\si)})$ is the smallest set of positive
integers that contains $\EP(\om)$ and satisfies assumptions (i), (ii),
(iii) of Lemma \ref{odo9}.
\end{proposition}

\begin{proof}
For any integer $n\ge1$ let $p_n$ be the smallest positive integer such
that $\om_n=\om_{n+kp_n}$ for $k=1,2,\dots$.  Clearly, $p_n\in\EP(\om)$.
Moreover, the sequence $p_1,p_2,\dots$ contains all essential partial
periods of $\om$.  For any $m\ge1$ let $q_m=\lcm(p_1,p_2,\dots,p_m)$.  Then
$\om_n=\om_{n+kq_m}$ for all $n\in\{1,2,\dots,m\}$ and $k\ge1$.  In
particular, the first $m$ elements of the sequence $\si^{q_m}\om$ are the
same as the first $m$ elements of $\om$.  It follows that
$\si^{q_m}\om\to\om$ as $m\to\infty$.

Suppose $n\in\CF(\si|_{Z(\om,\si)})$, i.e., the subshift
$\si|_{Z(\om,\si)}$ is a continuous extension of a cyclic permutation
$T:X\to X$ of order $n$.  Let $f:Z(\om,\si)\to X$ be a continuous map such
that $f\si=Tf$.  Since $\si^{q_m}\om\to\om$ as $m\to\infty$, we have that
$f(\si^{q_m}\om)\to f(\om)$ as $m\to\infty$.  But $X$ is a finite set with
the discrete topology so $f(\si^{q_m}\om)=T^{q_m}f(\om)$ actually coincides
with $f(\om)$ for large $m$.  It follows that $n$ divides $q_m$ for large
$m$.

Let $F_0$ denote the smallest set of positive integers that contains
$\EP(\om)$ and satisfies assumptions (i), (ii), (iii) of Lemma \ref{odo9}.
Clearly, a positive integer belongs to $F_0$ if and only if it divides
$\lcm(n_1,n_2,\dots,n_k)$ for some $n_1,n_2,\dots,n_k\in\EP(\om)$.  By the
above $\CF(\si|_{Z(\om,\si)})\subset F_0$.  On the other hand,
$\EP(\om)\subset\CF(\si|_{Z(\om,\si)})$ due to Lemma \ref{toe5}.  Then it
follows from Lemma \ref{odo8} that $F_0\subset\CF(\si|_{Z(\om,\si)})$.
Thus $\CF(\si|_{Z(\om,\si)})=F_0$.
\end{proof}

\begin{lemma}\label{toe7}
Let $\om$ be a Toeplitz sequence, $T:X\to X$ be a maximal odometer factor
of $\si|_{Z(\om,\si)}$, and $f:Z(\om,\si)\to X$ be a continuous map such
that $Tf=f\si$ on $Z(\om,\si)$.  Then, given $\om'\in Z(\om,\si)$, the
equation $f(\eta)=f(\om')$ has a solution $\eta$ different from $\om'$ if
and only if $\om'$ is not a Toeplitz sequence.
\end{lemma}

\begin{proof}
We have $X=\bZ_{n_1}\times\bZ_{n_2}\times\dots$ for some positive integers
$n_1,n_2,\dots$.  For any integer $k\ge1$ let $m_k=n_1n_2\cdots n_k$.
Given a sequence $x=(x_1,x_2,\dots)\in X$ and an integer $N\ge1$, the first
$k$ elements of the sequence $T^Nx$ match the first $k$ elements of $x$ if
and only if $m_k$ divides $N$.  Consequently, $T^{N_i}x\to x$ as
$i\to\infty$ if and only if each $m_k$ divides all but finitely many of the
numbers $N_1,N_2,\dots$.

First assume that $\om'=(\om'_1,\om'_2,\dots)\in Z(\om,\si)$ is not a
Toeplitz sequence.  Then there exists an index $m\ge1$ such that for any
integer $p\ge1$ a subsequence $\om'_m,\om'_{m+p},\om'_{m+2p},\dots$ is not
constant.  That is, $\om'_m\ne\om'_{m+N(p)p}$ for some $N(p)\ge1$.  Take
any limit point $\om''$ of a sequence $\si^{N(m_1)m_1}\om',
\si^{N(m_2)m_2}\om',\dots$.  By the above, $f(\si^{N(m_k)m_k}\om')=
T^{N(m_k)m_k}f(\om')\to f(\om')$ as $k\to\infty$.  Therefore
$f(\om'')=f(\om')$.  By construction, the sequences $\om''$ and $\om'$
differ at least in the $m$th element.

Now consider the case when $\om'$ is a Toeplitz sequence.  Take any
$\eta\in Z(\om,\si)$ such that $f(\eta)=f(\om')$.  According to Lemma
\ref{toe2}, the subshift $\si|_{Z(\om,\si}$ is minimal.  Hence
$Z(\om',\si)=Z(\om,\si)$.  Then for some integers $0\le N_1\le N_2\le\dots$
we have $\si^{N_i}\om'\to\eta$ as $i\to\infty$, which implies that
$T^{N_i}f(\om')=f(\si^{N_i}\om')\to f(\eta)=f(\om')$ as $i\to\infty$.  By
the above each $m_k$ divides all but finitely many of the numbers
$N_1,N_2,\dots$.  Since $\om'$ is a Toeplitz sequence, for any index
$m\ge1$ there exists $p\in\EP(\om')$ such that $\om'_m=\om'_{m+np}$,
$n=1,2,\dots$.  By Lemma \ref{toe5}, a cyclic permutation of order $p$ is a
continuous factor of the subshift $\si|_{Z(\om,\si)}$.  Then Lemma
\ref{odo6} implies that $p$ divides some $m_k$.  Hence $p$ divides $N_i$
for all sufficiently large $i$, which implies that $\om'_{m+N_i}=\om'_m$
for all sufficiently large $i$.  It follows that the $m$th element of the
sequence $\eta$ equals $\om'_m$.  As the choice of $m$ was arbitrary, we
conclude that $\eta=\om'$.
\end{proof}

Using notation of Lemma \ref{toe7}, let $\Omega_0$ be the set of all
Toeplitz sequences in $Z(\om,\si)$.  Clearly, $\Omega_0$ is invariant under
the shift.  According to Lemma \ref{toe7}, $\Omega_0$ is the largest subset
of $Z(\om,\si)$ such that the restriction of the mapping $f$ to $\Omega_0$
is one-to-one.  By Lemma \ref{toe4}, $\Omega_0$ is a residual subset of
$Z(\om,\si)$.  One might say that, generically, the dynamics of the
subshift $\si|_{Z(\om,\si)}$ is that of the odometer $T$.  However the
dynamics of this subshift as a whole can be much more complicated.  In the
next section we will consider a class of Toeplitz sequences $\om$ for which
the subshift $\si|_{Z(\om,\si)}$ is as close to an odometer as it gets.

\section{Substitution subshifts}\label{sub}

Given a finite alphabet $\cA$, we denote by $\cA^*$ the set of all finite
words in this alphabet.  Any $w\in\cA^*$ is simply an arbitrary finite
sequence of elements from $\cA$.  $\cA^*$ is a monoid with respect to
concatenation (the unit element is the empty word).  The concatenation of
a word $w\in\cA^*$ with a sequence $\om\in\cA^{\bN}$ is also naturally
defined, which gives rise to an action of the monoid $\cA^*$ on $\cA^{\bN}$.

Consider an arbitrary map $\tau:\cA\to\cA^*$.  Since $\cA^*$ is the free
monoid generated by $\cA$, the map $\tau$ can be extended to a homomorphism
of $\cA^*$ to itself.  Given a word $w\in\cA^*$, the word $\tau(w)$ is
obtained by substituting the word $\tau(a)$ for each occurence of any
letter $a\in\cA$ in $w$.  The latter procedure applies to infinite words as
well, which gives rise to a transformation of $\cA^{\bN}$ called a {\em
substitution}.  For convenience, we use the notation $\tau$ for both the
homomorphism of $\cA^*$ and the transformation of $\cA^{\bN}$.  Formally,
the map $\tau:\cA\to\cA^*$ is uniquely extended to a transformation of the
set $\cA^*\cup\cA^{\bN}$ such that $\tau(uw)=\tau(u)\tau(w)$ and
$\tau(u\om)=\tau(u)\tau(\om)$ for all $u,w\in\cA^*$ and $\om\in\cA^{\bN}$.

Let $\tau$ be a substitution on $\cA^{\bN}$ induced by a map $\tau:
\cA\to\cA^*$ as described above.  We assume $\tau$ to be non-degenerate in
that the word $\tau(b)$ is nonempty for any letter $b\in\cA$.  Suppose that
for some letter $a\in\cA$ the word $\tau(a)$ begins with $a$ and contains
more than one letter.  Consider the words $a,\tau(a),\tau^2(a),\tau^3(a),
\dots$.  It easily follows by induction that each word in this sequence is
a beginning of the next one and that each word is shorter than the next
one.  Therefore the finite words $a,\tau(a),\tau^2(a),\tau^3(a),\dots$ in a
sense converge to an infinite sequence $\om\in\cA^{\bN}$.  Namely, $\om$ is
the unique infinite sequence such that each $\tau^k(a)$ is a beginning of
$\om$.  By construction, the sequence $\om$ is a fixed point of the
substitution $\tau$.  The associated subshift $\si|_{Z(\om,\si)}$ is called
a {\em substitution subshift}.

In this paper, we are mostly interested in a substitution $\tau_G$ over the
alphabet $\cA=\{a,b,c,d\}$ that arises in the study of the Grigorchuk
group.  $\tau_G$ is defined by relations
$$
\tau_G(a)=aca,\quad \tau_G(b)=d,\quad \tau_G(c)=b,\quad \tau_G(d)=c.
$$
This substitution has a unique invariant sequence $\om^{(G)}$ that is the
limit of finite words $\tau_G^k(a)$, $k=1,2,\dots$.  The following lemma
provides a complete description of the sequence $\om^{(G)}=
\bigl(\om^{(G)}_1,\om^{(G)}_2,\dots\bigr)$.

\begin{lemma}\label{sub1}
Let $m=2^k(2n+1)$, where $k,n\ge0$ are integers.  If $k=0$ then
$\om^{(G)}_m=a$.  In the case $k>0$, we have $\om^{(G)}_m=d$, $c$, or $b$ if
the remainder of $k$ under division by $3$ is $0$, $1$, or $2$, respectively.
\end{lemma}

\begin{proof}
The first two letters of the infinite word $\om^{(G)}$ are $ac$.  For any
letter $l\in\{b,c,d\}$ we have $\tau_G(al)=acal'$, where $l'=\tau_G(l)$ is
another letter from $\{b,c,d\}$.  Using inductive argument, we derive from
these simple observations that $\om^{(G)}_m=a$ if and only if $m$ is odd.
Furthermore, $\om^{(G)}_{4n-2}=c$ and $\om^{(G)}_{4n}=
\tau_G\bigl(\om^{(G)}_{2n}\bigr)$ for $n=1,2,\dots$.

Let $m=2^k(2n+1)$, where $k,n\ge0$ are integers.  If $k=0$ then
$\om^{(G)}_m=a$ since $m$ is odd.  If $k=1$ then $\om^{(G)}_m=c$.  If $k>1$
then $\om^{(G)}_m=\tau_G^{k-1}\bigl(\om^{(G)}_{2(2n+1)}\bigr)=
\tau_G^{k-1}(c)$.  It remains to notice that $\tau_G^{k-1}(c)=d$, $c$, or
$b$ if the remainder of $k$ under division by $3$ is $0$, $1$, or $2$,
respectively.
\end{proof}

\begin{lemma}\label{sub2}
$\om^{(G)}$ is a Toeplitz sequence.  The set $\EP\bigl(\om^{(G)}\bigr)$ of
its essential partial periods consists of all powers of $2$.
\end{lemma}

Since $\om^{(G)}$ is a Toeplitz sequence, Lemma \ref{toe3} applies to it.
In fact, in this particular case a much stronger statement holds.

\begin{lemma}\label{sub3}
Let $\om=(\om_1,\om_2,\dots)\in Z\bigl(\om^{(G)},\si\bigr)$.  Given
positive integers $p$ and $m$, if $\om_m=\om_{m+p}=\om_{m+2p}=\om_{m+3p}$,
then $\om_{m+np}=\om_m$ for all $n\ge1$.
\end{lemma}

\begin{proofof}{Lemmas \ref{sub2} and \ref{sub3}}
For any integer $m>0$ let $d(m)$ denote the largest nonnegative integer
such that $2^{d(m)}$ divides $m$.  Clearly,
$d(m_1m_2)=d(m_1)+d(m_2)$ and $d(m_1+m_2)\ge\min(d(m_1),d(m_2))$.
Moreover, if $d(m_1)\ne d(m_2)$ then $d(m_1+m_2)=\min(d(m_1),d(m_2))$.

According to Lemma \ref{sub1}, any element $\om^{(G)}_m$ of the sequence
$\om^{(G)}$ is uniquely determined by $d(m)$.  In particular,
$\om^{(G)}_m\ne\om^{(G)}_n$ whenever $|d(m)-d(n)|=1$.

Take any integers $p,m\ge1$.  If $d(p)>d(m)$ then $d(m+np)=d(m)$ for all
$n\ge1$, which implies that $\om^{(G)}_{m+np}=\om^{(G)}_m$ for all $n\ge1$.

In the case $d(p)=d(m)$, we have $d(m+p)>d(m)$.  If, in addition,
$d(m+p)>d(m)+1=d(2p)$ then $d(m+3p)=d(2p)=d(m)+1$.  Therefore
$d(m+p)-d(p)=1$ or $d(m+3p)-d(p)=1$.

In the case $d(p)<d(m)$, we have $d(m+p)=d(p)$.  If, in addition,
$d(m)>d(p)+1=d(2p)$ then $d(m+2p)=d(2p)=d(p)+1=d(m+p)+1$.  Therefore
$d(m)-d(m+p)=1$ or $d(m+2p)-d(m+p)=1$.

By the above the equalities $\om^{(G)}_m=\om^{(G)}_{m+p}=\om^{(G)}_{m+2p}=
\om^{(G)}_{m+3p}$ cannot hold simultaneously if $d(p)\le d(m)$.

Denote by $\Omega$ the set of all sequences $\om=(\om_1,\om_2,\dots)\in
\cA^{\bN}$ such that, given arbitrary integers $p,m>0$, we have
$\om_{m+np}=\om_m$ for all $n\ge1$ whenever this holds for $1\le n\le3$.
It is easy to observe that the set $\Omega$ is shift-invariant and closed.
We have just shown that $\om^{(G)}\in\Omega$.  It follows that
$Z\bigl(\om^{(G)},\si\bigr)\subset\Omega$, which is exactly what Lemma
\ref{sub3} states.

We proceed to the proof of Lemma \ref{sub2}.  Given an integer $m\ge1$, let
$p=2^{d(m)+1}$.  Then $d(p)=d(m)+1>d(m)$.  This implies
$\om^{(G)}_{m+np}=\om^{(G)}_m$ for all $n\ge1$.  On the other hand, take
any $1\le p'<p$.  Then $d(p')<d(p)$ or, equivalently, $d(p')\le d(m)$.  By
the above $\om^{(G)}_{m+np'}\ne\om^{(G)}_m$ for some $1\le n\le 3$.
Since $m$ can be chosen arbitrarily, we conclude that $\om^{(G)}$ is a
Toeplitz sequence and $\EP\bigl(\om^{(G)}\bigr)=\{2^{d(m)+1}\mid m\ge1\}=
\{2^k\mid k\ge1\}$.
\end{proofof}

The remainder of this section is devoted to the study of the subshift
$\si|_{Z(\om,\si)}$ associated to a Toeplitz sequence $\om$ with
$\EP(\om)=\{2^k\mid k\ge1\}$.  We shall see that the latter condition
completely determines the dynamics of the subshift.  Note that, in general,
$\si|_{Z(\om,\si)}$ is not a substitution subshift.

\begin{proposition}\label{sub4}
Let $\om$ be a Toeplitz sequence such that $\EP(\om)=\{2^k\mid k\ge1\}$.
Then
\begin{itemize}
\item
the preimage under the subshift $\si|_{Z(\om,\si)}$ of any sequence
$\om'\in Z(\om,\si)$ is nonempty;
\item
there exists a unique sequence $\om^*\in Z(\om,\si)$ for which this
preimage contains more than one sequence;
\item
a sequence $\om'\in Z(\om,\si)$ is not a Toeplitz sequence if and only if
$\si^n\om'=\om^*$ for some $n\ge1$.
\end{itemize}
In the case $\om=\om^{(G)}$, we have $\om^*=\om^{(G)}$.  The preimage of
$\om^{(G)}$ under the subshift consists of three sequences $b\om^{(G)}$,
$c\om^{(G)}$, and $d\om^{(G)}$.
\end{proposition}

\begin{proof}
Let $\om=(\om_1,\om_2,\dots)$.  Take an arbitrary integer $k\ge1$.  Since
$2^k\in\EP(\om)$, there exist a letter $l_k\in\cA$ and an index $m_k\ge1$
such that $\om_{m_k+2^kn}=l_k$, $n=0,1,2,\dots$, but for any $1\le p<2^k$
the sequence $\om_{m_k},\om_{m_k+p},\om_{m_k+2p},\dots$ contains an element
different from $l_k$.  Lemma \ref{toe1} implies that $\om_n=l_k$ for every
positive $n\equiv m_k \bmod{2^k}$.  It follows that $m_K\not\equiv m_k
\bmod{2^k}$ if $K>k$.  As a consequence, the congruence classes
$m_1+2\bZ,m_2+2^2\bZ,\dots$ are pairwise disjoint.

Let $U_0=\bZ$ and $U_k=U_{k-1}\setminus (m_k+2^k\bZ)$ for $k\ge1$.  Since
$m_1+2\bZ,m_2+2^2\bZ,\dots$ are pairwise disjoint sets, it follows that
each $U_k$ is a congruence class modulo $2^k$.  That is, $U_k=M_k+2^k\bZ$
for a unique $1\le M_k\le 2^k$.  By construction, $M_k\equiv M_{k+1}
\bmod{2^k}$.  Hence $M_{k+1}=M_k$ or $M_{k+1}=M_k+2^k$.  In particular, the
sequence $M_1,M_2,\dots$ is nondecreasing.  Note that $\om_{m+2^k}=\om_m$
for any positive integer $m\not\equiv M_k\bmod{2^k}$.  This implies that
the sequence $\om_{M_k},\om_{M_k+2^k},\dots,\om_{M_k+2^kn},\dots$ contains
an element different from $\om_{M_k}$ as otherwise $\om$ will be periodic
(with period $2^k$), which is impossible since $\EP(\om)$ is an infinite
set.  On the other hand, since all essential partial periods of $\om$ are
powers of $2$, we do have $\om_{M_k+2^Kn}=\om_{M_k}$ for some $K>k$ and
$n=1,2,\dots$.  Then $M_k<M_K$.  Hence the sequence $M_1,M_2,\dots$ tends
to infinity.  As a consequence, the intersection of sets
$U_0,U_1,U_2,\dots$ contains no positive integer (and at most one
nonpositive integer).

The sequences $l_1,l_2,\dots$ and $M_1,M_2,\dots$ can be used to
reconstruct the sequence $\om$.  Namely, $\om_m=l_k$ whenever
$m\in U_{k-1}\setminus U_k$, that is, $m\not\equiv M_k\bmod{2^k}$ while
$m\equiv M_{k-1}\bmod{2^{k-1}}$.  Observe that the sequence of letters
$l_1,l_2,\dots$ is not eventually constant as otherwise $\om$ would be
periodic.

Now consider an arbitrary sequence $\om'=(\om'_1,\om'_2,\dots)\in
Z(\om,\si)$.  We have $\si^{n_i}\om\to\om'$ as $i\to\infty$ for some
$0\le n_1\le n_2\le\dots$.  Since $2^k\in\EP(\om)$ for any $k\ge1$, it
follows from Lemma \ref{toe5} that a cyclic permutation of order $2^k$ is a
continuous factor of the subshift $\si|_{Z(\om,\si)}$.  As a consequence,
the sequences $\si^{n_1}\om,\si^{n_2}\om,\dots$ may converge in
$Z(\om,\si)$ only if the numbers $n_1,n_2,\dots$ eventually have the same
remainder $r_k$ under division by $2^k$.  Furthermore, the remainder $r_k$
depends only on the limit sequence $\om'$ and not on the choice of $n_i$.
Therefore we have a unique, well defined integer $1\le M_k(\om')\le 2^k$
such that $M_k(\om')+r_k\equiv M_k\bmod{2^k}$.  Observe that $M_k(\om)=M_k$
and $M_k(\si\om')+1\equiv M_k(\om')\bmod{2^k}$.  By construction,
$r_{k+1}\equiv r_k\bmod{2^k}$, which implies that $M_{k+1}(\om')\equiv
M_k(\om')\bmod{2^k}$.  In other words, the congruence classes
$M_k(\om')+2^k\bZ$, $k\ge1$ are nested just as in the case $\om'=\om$.
There is more similarity.  For example, we know that a sequence
$\om_m,\om_{m+2^k},\dots,\om_{m+2^kn},\dots$ is constant whenever
$m\not\equiv M_k\bmod{2^k}$, while it is not constant for $m=M_k$.  It
follows that a sequence $\om'_m,\om'_{m+2^k},\dots,\om'_{m+2^kn},\dots$ is
constant whenever $m\not\equiv M_k(\om')\bmod{2^k}$.  On the other hand,
this sequence is not constant for $m=M_k(\om')$ as otherwise $\om'$ will be
periodic (with period $2^k$), which is impossible since the orbit
$\om',\si\om',\si^2\om',\dots$ is dense in the infinite set $Z(\om,\si)$
due to Lemma \ref{toe2}.  Actually, the only difference of the general case
from the case $\om'=\om$ is that the sequence $M_1(\om'),M_2(\om'),\dots$
may not tend to infinity.

The sequences $l_1,l_2,\dots$ and $M_1(\om'),M_2(\om'),\dots$ can be used
to reconstruct the sequence $\om'$.  We know that $\om_m=l_k$ whenever
$m\not\equiv M_k\bmod{2^k}$ and $m\equiv M_{k-1}\bmod{2^{k-1}}$.  It
follows that $\om'_m=l_k$ whenever $m\not\equiv M_k(\om')\bmod{2^k}$ and
$m\equiv M_{k-1}(\om')\bmod{2^{k-1}}$.  If $M_k(\om')\to\infty$ as
$k\to\infty$, then this is enough to reconstruct $\om'$.  Also, this
implies that $\om'$ is a Toeplitz sequence.  Otherwise, when the sequence
$M_1(\om'),M_2(\om'),\dots$ is eventually constant, we still need to
determine an element $\om_{M_\infty}$, where $M_\infty$ is the limit of
$M_k(\om')$ as $k\to\infty$.  By the way, $M_\infty$ completely determines
the sequence $M_1(\om'),M_2(\om'),\dots$ as $M_k(\om')\equiv M_\infty
\bmod{2^k}$, $k\ge1$.  Let $l_\infty(\om')=\om_{M_\infty}$.  Along with the
letter $l_\infty(\om')$, the sequences $l_1,l_2,\dots$ and
$M_1(\om'),M_2(\om'),\dots$ determine $\om'$ uniquely.  Notice that the
letter $l_\infty(\om')$ must occur infinitely often in the sequence
$l_1,l_2,\dots$.  Indeed, assume that some $a\in\cA$ never occurs in a
subsequence $l_k,l_{k+1},\dots$.  Then $\om_m=a$ implies
$\om_{m+2^{k-1}}=a$ for any $m\ge1$.  Consequently, $\om'_m=a$ implies
$\om'_{m+2^{k-1}}=a$ for any $m\ge1$.  It follows that
$\om'_{M_K(\om')}\ne a$ for all $K\ge k-1$.

Next consider an arbitrary sequence of integers $M'_1,M'_2,\dots$ such that
$1\le M'_k\le 2^k$ and $M'_{k+1}=M'_k\bmod{2^k}$ for all $k\ge1$.  We
associate to it another sequence of positive integers $R_1,R_2,\dots$
defined by $R_k=M_k-M'_k+2^k$.  For any $k\ge1$ the numbers
$R_k,R_{k+1},\dots$ have the same remainder $r_k$ under division by $2^k$.
Clearly, $M'_k+r_k\equiv M_k\bmod{2^k}$.  It follows from the above that
any limit point $\om'$ of the sequence $\si^{R_1}\om,\si^{R_2}\om,\dots$
satisfies $M_k(\om')=M'_k$ for all $k\ge1$.  In particular, there exists a
sequence $\om^*\in Z(\om,\si)$ such that $M_k(\om^*)=2^k$, $k=1,2,\dots$.
Such a sequence is unique since $2^k\to\infty$ as $k\to\infty$.

Further, assume that the sequence $M'_1,M'_2,\dots$ from the previous
paragraph is eventually constant and denote by $M_\infty$ its limit.  By
the above there is at least one $\om'=(\om'_1,\om'_2,\dots)\in Z(\om,\si)$
such that $M_k(\om')=M'_k$ for all $k\ge1$.  Pick any letter $l\in\cA$ that
occurs infinitely often in the sequence $l_1,l_2,\dots$.  Let $\om''$ be
the sequence obtained from $\om'$ by substituting $l$ for a single element
$\om'_{M_\infty}$.  We shall show that $\om''\in Z(\om,\si)$.  Given an
integer $k\ge1$, a sequence $\om'_m,\om'_{m+2^k},\dots,\om'_{m+2^kn},\dots$
is constant whenever $m\not\equiv M_k(\om')\equiv M_\infty\bmod{2^k}$.  On
the other hand, in the case $m=M_\infty$ this sequence is not constant as
it contains all letters that occur in the sequence $l_{k+1},l_{k+2},\dots$.
In particular, $\om'_{M_\infty+2^kn}=l$ for some $n\ge0$.  Then the first
$2^k$ elements of the shifted sequence $\si^{2^kn}\om'$ coincide with the
first $2^k$ elements of $\om''$.  Since $k$ can be chosen arbitrarily
large, it follows that $\om''\in Z(\om,\si)$.  Note that
$M_k(\om'')=M_k(\om')=M'_k$ for all $k\ge1$ and that $l_\infty(\om'')=l$.

We already know that a sequence $\om'\in Z(\om,\si)$ is a Toeplitz sequence
if $M_k(\om')\to\infty$ as $k\to\infty$.  Let us show that $\om'$ is not a
Toeplitz sequence if the sequence $M_1(\om'),M_2(\om'),\dots$ is
eventually constant.  As follows from the above, in this case there exists
a sequence $\om''\in Z(\om,\si)$ that differs from $\om'$ in a single
element.  Then $\si^n\om'=\si^n\om''$ for some $n\ge1$.  Now we are going
to apply Lemma \ref{toe7}.  Let $T:X\to X$ be a maximal odometer factor of
the subshift $\si|_{Z(\om,\si)}$ and $f:Z(\om,\si)\to X$ be a continuous
map such that $Tf=f\si$ on $Z(\om,\si)$.  Since $\si^n\om'=\si^n\om''$, we
obtain $T^n(f(\om'))=T^n(f(\om''))$.  But every odometer is a one-to-one
mapping, which implies that $f(\om')=f(\om'')$.  By Lemma \ref{toe7},
neither $\om'$ nor $\om''$ is a Toeplitz sequence.

Recall that $M_k(\si\om')+1\equiv M_k(\om')\bmod{2^k}$ for all $k\ge1$ and
$\om'\in Z(\om,\si)$.  Consequently, for any $n\ge1$ we have
$M_k(\si^n\om')+n\equiv M_k(\om')\bmod{2^k}$.  It follows that the sequence
$M_1(\om'),M_2(\om'),\dots$ has a finite limit $n$ if and only if
$M_k(\si^n\om')=2^k$ for all $k\ge1$.  An equivalent condition is that
$\si^n\om'=\om^*$.  Thus a sequence $\om'\in Z(\om,\si)$ is not a Toeplitz
sequence if and only if $\si^n\om'=\om^*$ for some $n\ge1$.

By the above two sequences $\om',\om''\in Z(\om,\si)$ coincide if and only
if $M_k(\om')=M_k(\om'')$ for $k=1,2,\dots$ and either the sequence
$M_1(\om'),M_2(\om'),\dots$ tends to infinity, or else it is eventually
constant and $l_\infty(\om')=l_\infty(\om'')$.  Take any $\om'\in
Z(\om,\si)$ and consider a unique sequence of integers $M'_1,M'_2,\dots$
such that $1\le M'_k\le 2^k$ and $M_k(\om')+1\equiv M'_k\bmod{2^k}$ for all
$k\ge1$.  The congruency classes $M'_k+2^k\bZ$, $k\ge1$ are nested since
the congruency classes $M_k(\om')+2^k\bZ$, $k\ge1$ are nested.  As shown
above, this implies the existence of an $\om''\in Z(\om,\si)$ such that
$M_k(\om'')=M'_k$, $k\ge1$.  Let $S$ be the set of all such $\om''$.  The
set $S$ contains a single point if $M'_k\to\infty$ as $k\to\infty$.
Otherwise $S$ contains several points; they are distinguished by
$l_\infty(\om'')$, which can be any letter that occurs infinitely often in
the sequence $l_1,l_2,\dots$.  Clearly, a sequence $\om''\in Z(\om,\si)$
belongs to $S$ if and only if $M_k(\si\om'')=M_k(\om')$ for all $k\ge1$.
If $M_k(\om')\to\infty$ as $k\to\infty$, then $S$ is exactly the preimage
of $\om'$ under the subshift $\si|_{Z(\om,\si)}$.  If, in addition,
$\om'\ne\om^*$, then $M'_1,M'_2,\dots$ tend to infinity as well so that the
set $S$ consists of a single point.  In the case $\om'=\om^*$, all $M'_k$
are equal to $1$ so that $S$ contains more than one point.  Finally,
consider the case when the sequence $M_1(\om'),M_2(\om'),\dots$ is
eventually constant.  Then the sequence $M'_1,M'_2,\dots$ is also
eventually constant.  For any $\om''\in S$ we have
$l_\infty(\si\om'')=l_\infty(\om'')$.  Hence the preimage of $\om'$ under
the subshift consists of a unique $\om''\in S$ such that
$l_\infty(\om'')=l_\infty(\om')$.

The general part of the proposition is proved.  It remains to consider a
particular case $\om=\om^{(G)}$.  According to Lemma \ref{sub1}, in this
case the sequence $l_1,l_2,\dots$ is eventually periodic with period $3$:
$a,c,b,d,c,b,d,\dots$.  Moreover, $M_k=2^k$ for any $k\ge1$.  This means
that $\om^{(G)}=\om^*$.  By the above the preimage of $\om^{(G)}$ under the
subshift $\si|_{Z(\om^{(G)},\si)}$ consists of all sequences of the form
$l\om^{(G)}$, where $l\in\cA$ is a letter that occurs infinitely often in
the sequence $l_1,l_2,\dots$.  We have three such letters: $b$, $c$, and
$d$.
\end{proof}

\begin{theorem}\label{sub5}
Let $\om$ be a Toeplitz sequence such that $\EP(\om)=\{2^k\mid k\ge1\}$.
Then the maximal odometer factor of the subshift $\si|_{Z(\om,\si)}$ is the
binary odometer.

Moreover, the subshift $\si|_{Z(\om,\si)}$ is, up to a countable set,
continuously conjugated to the binary odometer.  To be precise, there exist
a countable set $\Omega_1\subset Z(\om,\si)$ and a continuous mapping $f$
of $Z(\om,\si)$ onto the ring $\cZ_2$ of dyadic integers such that
$f(\si\om')=f(\om')+1$ for all $\om'\in Z(\om,\si)$, the complement
$Z(\om,\si)\setminus\Omega_1$ is shift invariant and $f$ is one-to-one when
restricted to $Z(\om,\si)\setminus\Omega_1$.
\end{theorem}

\begin{proof}
Proposition \ref{toe6} implies that $\CF\bigl(\si|_{Z(\om,\si)}\bigr)=
\{2^k\mid k\ge0\}$.  That is, a nontrivial cyclic permutation is a
continuous factor of the subshift $\si|_{Z(\om,\si)}$ if and only if its
order is a power of $2$.  Let $T_0$ denote the odometer on
$\bZ_2\times\bZ_2\times\dots$.  By Lemma \ref{odo6}, $\CF(T_0)=
\{2^k\mid k\ge0\}=\CF\bigl(\si|_{Z(\om,\si)}\bigr)$.  Then it follows from
Proposition \ref{odo11} that $T_0$ is a maximal odometer factor of
$\si|_{Z(\om,\si)}$.  Finally, $T_0$ is continuously conjugated to the
binary odometer due to Lemma \ref{odo12}.

Thus there exists a continuous onto mapping $f:Z(\om,\si)\to\cZ_2$ such
that $f(\si\om')=f(\om')+1$ for all $\om'\in Z(\om,\si)$.  Denote by
$\Omega_0$ the set of all Toeplitz sequences in $Z(\om,\si)$.  Obviously,
the set $\Omega_0$ is shift invariant.  According to Lemma \ref{toe7}, the
map $f$ is one-to-one when restricted to $\Omega_0$.  By Proposition
\ref{sub4}, there exists a sequence $\om^*\in Z(\om,\si)$ such that any
$\om'\in\Omega_1=Z(\om,\si)\setminus\Omega_0$ satisfies $\si^n\om'=\om^*$
for some $n\ge1$.  For any fixed $n$ there are only finitely many sequences
$\om'$ satisfying this relation.  Therefore $\Omega_1$ is a countable set.
\end{proof}

Theorem \ref{sub5} applies to the sequence $\om^{(G)}$.  Hence there exists
a continuous map $f_G:Z(\om^{(G)},\si)\to\cZ_2$ such that $f_G(\si\om')=
f_G(\om')+1$ for all $\om'\in Z(\om^{(G)},\si)$.  According to Lemma
\ref{odo3}, we can choose this map so that $f_G(\om^{(G)})=0$; then it is
uniquely determined.  Recall that any dyadic integer $z\in\cZ_2$ is
uniquely expanded into a series $\sum_{i=1}^\infty n_i2^{i-1}$, where
$n_i\in\{0,1\}$.  Therefore the map $f_G$ can be regarded as a symbolic map
that takes an infinite word over the alphabet $\cA$ and assigns to it an
infinite word over the alphabet $\{0,1\}$.  The proof of Proposition
\ref{sub4} suggests an algorithm for effective computation of $f_G$.  Take
an arbitrary sequence $\om'=(\om'_1,\om'_2,\dots)\in Z(\om^{(G)},\si)$.
For any integer $k\ge1$ there exists a unique integer $1\le M_k\le 2^k$
such that the subsequence $\om'_{M_k},\om'_{M_k+2^k},\dots,\om'_{M_k+2^kn},
\dots$ is not constant.  Note that $0\le 2^k-M_k<2^k$.  The dyadic
expansion of the number $2^k-M_k$ gives us the first $k$ coefficients in
the dyadic expansion of $f_G(\om')$.  That is, if $f_G(\om')=
\sum_{i=1}^\infty n_i2^{i-1}$, where $n_i\in\{0,1\}$, then $2^k-M_k=
\sum_{i=1}^k n_i2^{i-1}$.  By Lemma \ref{sub3}, $\om'_{M_k+2^kn}\ne
\om'_{M_k}$ already for some $1\le n\le 3$.  It follows that the first $k$
coefficients in the dyadic expansion of $f_G(\om')$ depend only on the
first $2^{k+2}$ elements of the sequence $\om'$.

\begin{theorem}\label{sub6}
Let $\om$ be a Toeplitz sequence such that $\EP(\om)=\{2^k\mid k\ge1\}$.
Then
\begin{itemize}
\item
there exists a unique Borel probability measure $\mu$ on $Z(\om,\si)$
invariant under the subshift $\tilde\si=\si|_{Z(\om,\si)}$.
\item
The subshift $\tilde\si$ is ergodic with respect to the measure $\mu$.
Moreover, every orbit of $\tilde\si$ is uniformly distributed with respect
to $\mu$.
\item
$\tilde\si$ has purely point spectrum, the eigenvalues being all roots of
unity of order $1,2,\dots,2^k,\dots$.  Each eigenvalue is simple.
\item
All eigenfunctions of $\tilde\si$ are continuous.
\end{itemize}
\end{theorem}

\begin{proof}
Let $\Omega_1$ be the smallest subset of $Z(\om,\si)$ that contains all
non-Toeplitz sequences in $Z(\om,\si)$ and is both forward and backward
invariant under the subshift $\tilde\si$.  Proposition \ref{sub4} implies
that the set $\Omega_1$ is countable.  Denote by $\Omega_0$ the complement
of $\Omega_1$ in $Z(\om,\si)$.

By Theorem \ref{sub5}, the binary odometer is a maximal odometer factor of
$\tilde\si$.  Let $f:Z(\om,\si)\to\cZ_2$ be a continuous function such that
$f(\si\om')=f(\om')+1$ for all $\om'\in Z(\om,\si)$.  Lemma \ref{toe7}
implies that $f$ is one-to-one when restricted to $\Omega_0$.  Let
$X_0=f(\Omega_0)$ and denote by $h$ the inverse of the restriction of $f$
to $\Omega_0$.  By construction, the countable set $X_1=f(\Omega_1)$ is the
complement of $X_0$ in $\cZ_2$.  It follows that the mapping
$h:X_0\to\Omega_0$ is continuous.  By Proposition \ref{eigen2}, there
exists a unique Borel probability measure $\mu_0$ on $\cZ_2$ that is
invariant under the binary odometer.  For any Borel set $B\subset
Z(\om,\si)$ let $\mu(B)=\mu_0(h^{-1}(B\cap\Omega_0))$.  Then $\mu$ is
a Borel measure on $Z(\om,\si)$ invariant under the subshift $\tilde\si$.
Note that $\mu_0(X_1)=0$ since $X_1$ is a countable set and the binary
odometer has no finite orbits.  Therefore $\mu$ is a probability measure.

Suppose that $\mu'$ is a Borel probability measure on $Z(\om,\si)$
invariant under the subshift $\tilde\si$.  For any Borel set
$b\subset\cZ_2$ let $\mu'_0(b)=\mu'(f^{-1}(b))$.  It is easy to check that
$\mu'_0$ is a Borel probability measure on $\cZ_2$ invariant under the
odometer.  Hence $\mu'_0=\mu_0$.  It follows that $\mu'$ coincides with
$\mu$ on the set $\Omega_0$.  Since $\Omega_1$ is a countable set and the
subshift $\tilde\si$ has no finite orbits, we obtain $\mu'(\Omega_1)=
\mu(\Omega_1)=0$ so that $\mu'=\mu$.
Thus $\mu$ is the only shift-invariant Borel probability measure on
$Z(\om,\si)$.  It follows that the subshift $\tilde\si$ is ergodic with
respect to $\mu$, moreover, each orbit of $\tilde\si$ is uniformly
distributed in $Z(\om,\si)$.

Since $\mu(\Omega_1)=\mu_0(X_1)=0$, the subshift $\tilde\si$ and the binary
odometer $T_0$ are isomorphic as measure-preserving transformations.  Then
Proposition \ref{eigen2} implies that $\tilde\si$ has pure point spectrum,
the eigenvalues of $\tilde\si$ are all $n$th roots of unity, where $n$ runs
through $\CF(T_0)$, and all eigenvalues are simple.  The set $\CF(T_0)$
coincides with $\CF(\tilde\si)$ due to Proposition \ref{odo11}.  Further,
$\CF(\tilde\si)=\{2^k\mid k\ge0\}$ due to Proposition \ref{toe6}.

To prove that all eigenfunctions of $\tilde\si$ are continuous, it is
enough to show that for any root of unity $\zeta$ of order $2^k$ there
exists an associated continuous eigenfuction of $\tilde\si$.  Let
$\zeta_0=\exp(2\pi i/2^k)$.  Then $\zeta=\zeta_0^m$ for some $m\ge1$.  By
Lemma \ref{eigen1}, there is a continuous eigenfunction $\phi$ of
$\tilde\si$ associated with the eigenvalue $\zeta_0$.  Then $\phi^m$ is
also a continuous eigenfunction and its eigenvalue is $\zeta_0^m=\zeta$.
\end{proof}

\begin{proofof}{Theorems \ref{main1} and \ref{main2}}
According to Lemma \ref{sub2}, the sequence $\om^{(G)}$ is a Toeplitz
sequence with $\EP(\om^{(G)})=\{2^k\mid k\ge1\}$.  Therefore Theorem
\ref{main1} is a particular case of Theorem \ref{sub5} while Theorem
\ref{main2} is a particular case of Theorem \ref{sub6}.
\end{proofof}

\bigskip

{\sc
\begin{raggedright}
Department of Mathematics\\
Texas A\&M University\\
College Station, TX 77843--3368
\end{raggedright}
}

\medskip\noindent
{\it E-mail address:} {\tt yvorobet@math.tamu.edu}

\end{document}